\documentclass{amsart}
\usepackage{tabu}
\usepackage[margin=95pt]{geometry}
\usepackage{enumitem}
\allowdisplaybreaks

\usepackage{amssymb}
\usepackage{amsmath}
\usepackage{amscd}
\usepackage{amsbsy}
\usepackage{comment}
\usepackage{enumitem}
\usepackage[matrix,arrow]{xy}
\usepackage{hyperref}
\newcommand{\Q}{{\mathbb Q}}
\newcommand{\Z}{{\mathbb Z}}

\newcommand{\R}{{\mathbb R}}

\def\mod#1{{\ifmmode\text{\rm\ (mod~$#1$)}
\else\discretionary{}{}{\hbox{ }}\rm(mod~$#1$)\fi}}

\begin {document}

\newtheorem{thm}{Theorem}

\theoremstyle{definition}

\theoremstyle{remark}

\title[Pencils of norm-form equations]
{A note on pencils of norm-form equations}

\author{Prajeet Bajpai}
\address{Department of Mathematics, University of British Columbia, Vancouver, B.C., V6T 1Z2 Canada}
\email{prajeet@math.ubc.ca}

\author[Michael Bennett]{Michael A. Bennett}
\address{Department of Mathematics, University of British Columbia, Vancouver, B.C., V6T 1Z2 Canada}
\email{bennett@math.ubc.ca}

\thanks{The second author was supported by NSERC}

\date{\today}


\begin {abstract}
We find all solutions to the parametrized family of norm-form equations
$$
x^3-(t^3-1)y^3+3(t^3-1)xy+(t^3-1)^2 = \pm 1,
$$
studied by Amoroso, Masser and Zannier. Our proof relies upon an appeal to lower bounds for linear forms in logarithms and various elementary arguments.
\end {abstract}
\maketitle

\section{Introduction}

If $d \geq 3$ is an integer and $A_i(t) \in \mathbb{Z}[t]$, for $1 \leq i \leq d-1$, consider the binary form
\[
\Phi_t(x,y) = x(x-A_1(t) y )\cdots( x-A_{d-1}(t) y ) + y^d 
\]
and associated  Thue equations 
\begin{equation} \label{Thue}
\Phi_t(x,y) = 1 
\end{equation}
to be solved in integers $x,y,t$. In 1993, Thomas \cite{Th93} conjectured that if the degrees of the polynomials $A_i(t)$ are distinct, then this family of equations (parameterized by $t$) has only the ``trivial'' solutions
\[
(x,y) = (1,0), (0,1), (A_1(t),1), \ldots  (A_{d-1}(t),1),
\]
provided that  $t$ is sufficiently large. These are specializations of solutions to $\Phi_t(x,y) = 1$ with $x,y \in \Z[t]$, so Thomas's conjecture asserts that all integer solutions (for large $t$) come from such ``functional'' solutions, and moreover that the only functional solutions are the trivial ones. In 2007, Ziegler \cite{Zie} found a counterexample to this second aspect of Thomas's conjecture in the case $d=3$ and $A_1(t) = t, A_2(t) = t^4 + 3t$. Under these assumptions, we find an additional functional solution with
\[
x = t^9 + 3t^6 + 4t^3 + 1  \; \mbox{ and } \;  y = t^8 + 3t^5 + 3t^2
\]
to the equation $x(x-A_1(t)y)(x-A_2(t) y) + y^3 = 1$. For a sampling of the extensive  literature related to solving parametrized families of  Thue equations, the reader might wish to consult \cite{BeGh}, \cite{Heu}, \cite{Heu2}, \cite{HTZ}, \cite{Mig3}, \cite{MiT}, \cite{Tho}, \cite{Th93} and \cite{Wa5}. Fundamentally underpinning all of these results are lower bounds for linear forms in complex logarithms. That such techniques permit the effective solution of equations like (\ref{Thue}) for fixed $t$ is well known -- the challenge for such families is to demonstrate that, for suitably large $t$, no additional solutions accrue.

Recently, Amoroso, Masser and Zannier \cite{AMZ2} have shown that, at least for $d=3$, Thomas's equation admits only functional solutions when $t$ is larger than some effective $t_0$. Combined with work of  Ziegler \cite{Zie}, this confirms Thomas's conjecture when $d=3$ and $\deg(A_2) > 34\deg(A_1)$.
The techniques of  \cite{AMZ2}, based upon height bounds for solutions to equations in multiplicative tori, actually extend to quite general families of norm-form equations which are not of Thue type and {\it a priori} cannot be resolved through appeal to bounds for linear forms in logarithms.
As an example, in \cite{AMZ} and \cite{AMZ2}, Amoroso, Masser and Zannier consider the family of  inhomogeneous (non-Thue) norm-form equations
\[
x^3-(t^3-1)y^3+3(t^3-1)xy+(t^3-1)^2 = 1
\]
and prove that there exists an effective $t_0$ such that all integer solutions to the above equation with $t>t_0$ are of the form
\[
(x,y) = (t^2,2t),\quad (t^2,-t).
\]
Their proof, relying  upon the main theorem from \cite{AMZ},  does not provide any information on solutions to these equations for $t \le t_0$. In this short note, we prove that the functional solutions are in fact the only solutions to a slight generalization of  the above equation, for \emph{all} values of $t$. 
\begin{thm} \label{thm-main}
Let $x, y$ and $t$ be integers with $t>1$. If 
\begin{equation}\label{cubicnorm}
x^3-(t^3-1)y^3+3(t^3-1)xy+(t^3-1)^2 = \pm 1,
\end{equation}
then $(x,y) = (t^2,2t)$ or $(x,y) = (t^2,-t)$.
\end{thm}

As far as we are aware, this is the first example of a parametrized family of non-Thue norm-form equations to be solved completely.
Our proof will in fact rely upon bounds for linear forms in logarithms, together with 
an elementary argument to 
 ensure suitably rapid growth (which Thomas \cite{Th93} terms {\it stable growth}) of the exponent of the fundamental unit in the corresponding cubic field, coming from a solution to (\ref{cubicnorm}). The approach of this paper may be more generally applicable, including to the other families of norm-form equations considered in \cite{AMZ2}, but this is currently unclear.

\section{Proof of Theorem \ref{thm-main}}

Let $t>1$ be an integer,  define $s=\sqrt[3]{t^3-1}$ and 
set $\omega=e^{2\pi i/3}$. Equation \eqref{cubicnorm} then factorizes over $\Q(s,\omega)$ as
\[
(x-sy+s^2)(x-\omega sy+\omega^2 s^2)(x-\omega^2 sy+\omega s^2) = \pm 1,
\]
whence $x-sy+s^2$ is necessarily a unit in $\Q (s)$. For this latter (complex cubic) field, it is easy to directly determine a fundamental unit.  Since we have
\[
\mathrm{disc}(\Q(s)) = -27(t^3-1)^2,
\]
it follows from Artin's inequality (see e.g. \cite{Williams}, p. 370) that if $u>1$ is a unit in the ring of integers of this number field, then necessarily
\[
4u^3+27 > \left| \mathrm{disc}(\Q(s)) \right| = 27(t^3-1)^2.
\]
Now $t-s$ is a unit, and $(t-s)^{-1} = t^2+s^2+st > 1$. Further,
\[
4(t^2+s^2+st)^{3/2}+27 < 12\sqrt3 t^3 +27 < 27(t^3+1) < 27(t^3-1)^2
\]
for all integers $t>1$, and hence if $t^2+s^2+st=u^n$ for some unit $u > 1$ and integer $n \geq 2$, 
$$
27(t^3-1)^2 < 4u^3+27 = 4 (t^2+s^2+st)^{3/n}+27 \leq 4(t^2+s^2+st)^{3/2}+27 < 27(t^3-1)^2,
$$
a contradiction.
It follows that necessarily $n=1$ and so $t-s$ is a fundamental unit in $\Q(s)$. We may therefore write
\begin{equation} \label{start}
x-sy+s^2 = (-1)^\delta (t-s)^m
\end{equation}
for some integer $m$ and $\delta \in \{ 0, 1 \}$. Replacing $s$ by $\omega s$ and multiplying through by $\omega$, and then doing the same with $\omega$ replaced by $\omega^2$ implies that we also have
$$
\omega x- \omega^2 sy+s^2 = (-1)^\delta  (t-\omega s)^m \; \mbox{ and } \; \omega^2 x-\omega sy+s^2 = (-1)^\delta (t-\omega^2 s)^m.
$$
Summing these three equations, we find that
\begin{equation}\label{masser}
(t-s)^m + \omega(t-\omega s)^m + \omega^2(t-\omega^2 s)^m = (-1)^\delta 3s^2,
\end{equation}
a unit equation in the number field $\Q(s,\omega)$.

Generally, we cannot effectively solve a four-term unit equation like this, unless it satisfies certain extremely restrictive assumptions (see e.g. work of Vojta \cite{Voj}). 
Since, in our situation, $s$ alone does not generate a Galois field over $\Q$,  we cannot obtain a second independent unit equation to complement \eqref{masser} in a straightforward manner, whence the methods of \cite{Voj} do not directly apply.
Happily, though, as we shall see, equation (\ref{start}) has suitably pleasant arithmetic properties to enable us to solve (\ref{masser}) completely.

Let us begin by noting that the exponent $m=-1$ in (\ref{start})  corresponds to the functional solution $(x,y) = (t,-t)$ and $m=2$ corresponds to $(x,y)=(t,2t)$. If $m \in \{ 0, 1 \}$ then $s=0$ which contradicts $t>1$. Since
\[
|t-\omega s| =|t-\omega^2 s| =  \sqrt{t^2 + s^2 + st} > \sqrt{3s^2},
\]
then $m\le -2$ gives
\[
|\omega(t-\omega s)^m+\omega^2 (t-\omega^2 s)^m| \leq |t-\omega s|^m + |t-\omega^2 s|^m = 2 (t^2+s^2+st)^{m/2}
\le \frac{2}{3s^2}
\]
and
\[
|(t-s)^m| = |(t^2+s^2+st)^{-m}| > 9 s^{4}
\]
and so
\[
\big| (t-s)^m + \omega(t-\omega s)^m + \omega^2(t-\omega^2 s)^m \, \big| >  9s^{4} - \frac{2}{3s^2},
\]
which, since $t \geq 2$, contradicts (\ref{masser}) and again yields no solutions to equation (\ref{cubicnorm}). We will thus assume, here and henceforth, that $m\ge 3$.

\subsection{An upper bound on $m$}
Fix a determination of the complex logarithm such that $\log z = \log|z|+i\theta$ with $-\pi<\theta \le \pi$ and $z=|z|e^{i\theta}$. We note that, for $t\ge 2$, we have $|t-s|<1$. We will now appeal to lower bounds for  linear forms in complex logarithms to derive a contradiction for large $m$. We have
\[
\left| \omega(t-\omega s)^m + \omega^2(t-\omega^2 s)^m \right| = |t^2+s^2+st|^\frac{m}{2}\left| \frac{\omega(t-\omega s)^m}{\omega^2(t-\omega^2 s)^m} + 1 \right| 
\]
and noting that $(\omega(t-\omega s)^m)/(\omega^2(t-\omega^2 s)^m)$ lies on the unit circle,  we find that
\[
\frac12 \left| \log\left( \frac{-\omega(t-\omega s)^m}{\omega^2(t-\omega^2 s)^m} \right) \right|<  \left| \frac{-\omega(t-\omega s)^m}{\omega^2(t-\omega^2 s)^m} - 1 \right|.
\]
Set
\[
\Lambda_1 = \log\left( \frac{-1}\omega \left(\frac{t-\omega s}{t-\omega^2 s}\right)^m \right) = m\log\left( \frac{t-\omega s}{t-\omega^2 s} \right) - b\left(  \frac{2\pi i}{3}  \right)
\]
where we choose the principal branch of the logarithm and $b$ is taken so that $|\Lambda_1|$ is minimal. Since we have
\begin{equation} \label{upper}
\left| \Lambda_1 \right| < 2 \frac{\left| \omega(t-\omega s)^m + \omega^2(t-\omega^2 s)^m \right|}{ |t^2+s^2+st|^\frac{m}{2}}
\leq \frac{6s^2 + 2 (t-s)^m}{ |t^2+s^2+st|^\frac{m}{2}} < \frac{6}{3^{m/4}} s^{2-m} + \frac{2(t-s)^m}{\sqrt{3}s^m},
\end{equation}
and $m \geq 3$, the fact that $|\mathrm{Im}\log(t-\omega^2 s) | = |\mathrm{Im}\log(t-\omega s) |< \pi i/3$ implies that  necessarily  $|b|\le \frac m2+1$. Since $m>2$, in particular we have $|b| < m$. To derive a lower bound for $|\Lambda_1|$, we appeal to the following result, a special case of Corollary 1 of \cite{Lau}.
\begin{thm}[Laurent, 2008]\label{laurent}
Let  $\alpha_1$ and $\alpha_2$ be multiplicatively independent non-zero algebraic numbers with absolute logarithmic heights $h(\alpha_1), h(\alpha_2)$ respectively, and let $\log(\alpha_1), \log(\alpha_2)$ be any determination of their logarithms. Let $b_1$ and $b_2$ be positive integers, and set 
$$
D = [\Q(\alpha_1,\alpha_2) : \Q  ]/[\R=(\alpha_1,\alpha_2) : \R ]
$$
and
\[
b' = \frac{b_1}{D\log A_2} + \frac{b_2}{D\log A_1}
\]
where $A_1,A_2$ are real numbers $>1$ such that
\[
\log A_i \ge \max\{ h(\alpha_i), |\log\alpha_i|/D, 1/D  \} \quad (i = 1,2).
\]
Then the linear form
\[
\Lambda = b_2\log\alpha_2 - b_1\log\alpha_1
\]
satisfies the lower bound
\[
\log |\Lambda| \ge -22.8D^4(\max\{\log b' + 0.21, 30/D, 1\})^2\log A_1 \log A_2.
\]
\end{thm}

To apply this theorem in our case we first take
\[
D = [\Q(s,\omega):\Q]/[\R(s,\omega):\R] = 6/2 = 3
\]
and note that
\[
h\left( \frac{t-\omega s}{t-\omega^2s} \right) = \frac 13 \left( \log \left| \frac{t-\omega^2 s}{t-s} \right| \right) =  \frac12 \log (t^2+s^2+st),
\]
whence $ \frac{t-\omega s}{t-\omega^2s}$ and $e^{2\pi i/3}$ are multiplicatively independent and we can take
\begin{gather*}
\log A_1 = \max\left\{ h\left( \frac{t-\omega s}{t-\omega^2s} \right), \frac 13 \log\left|\frac{t-\omega s}{t-\omega^2s} \right|, \frac 13 \right\} = \frac12 \log (t^2+s^2+st)\\
\log A_2 = \max\left\{ h(\omega),  |\log \omega|/3, 1/3 \right\} = \frac{2\pi}9 .
\end{gather*}
This gives 
\[
22.8 D^4 \log A_1 \log A_2 < 644.66 \log(t^2+s^2+st).
\]
Next, we see
\[
 b' = \frac{2m}{3\log(t^2+s^2+st)} + \frac{3m}{2\pi} < 2m
\]
and so for  $m\ge 8928$ we have
\[
\max\{\log b' + 0.21, 10, 1\} < \log(2m) + 0.21.
\]
Therefore Theorem \ref{laurent} gives, assuming $m\ge 8928$, that
\[
\log|\Lambda_1| \ge -644.66 \, (\log (2m) + 0.21)^2 \log(t^2+s^2+st)
\]
and thus, from (\ref{upper}), we find that 
$$
\frac{m}{2} <  644.66 \, (\log (2m) + 0.21)^2 + \frac{\log (6s^2+2(t-s)^m)}{\log(t^2+s^2+st)} <  644.66 \, (\log (2m) + 0.21)^2 + 1.3,
$$
since $t \geq 2$ and $m \geq 3$.
We thus obtain an absolute upper bound
\begin{equation} \label{upperm}
m \le 225676,
\end{equation}
valid for all integers $t \geq 2$.

\subsection{Upper bound on $t$ and ruling out exceptional solutions}
Expanding $(t-s)^m$ as a binomial series and considering the coefficient of $s^2$ in equation \eqref{start}, we find that
\begin{equation} \label{expand}
\pm 1 = \sum_{i=0}^{[(m-2)/3]} \binom{m}{3i+2} t^{m-2-3i} (1-t^3)^i.
\end{equation}
A short argument modulo $3$ yields 
$$
\binom{m}{3i+2} \equiv 0 \mod{3}
$$
for all $m \equiv 0, 1 \mod{3}$ and so we necessarily have $m \equiv -1 \mod{3}$. The right-hand-side of (\ref{expand}) is thus congruent to $1$ modulo $t^3$, so that we may discard the ``$-$'' sign. Considering (\ref{expand}) modulo $t^3-1$, it follows that
\begin{equation}\label{smallhammer}
\binom{m}{2} t^{m-2} \equiv \binom{m}{2} \equiv  1 \mod{t^3-1}
\end{equation}
and so, from $m > 2$, that
$$
\binom{m}{2} \geq t^3 \; \Rightarrow \; m >  \sqrt{2} \,  t^{3/2}.
$$
Combining this with (\ref{upperm}) implies the inequality $t \le 2942$.

To sharpen this, we will improve the lower bound $m > \sqrt2 t^{3/2}$. Let us write $\binom{m}{2} =k(t^3-1)+1$, where $k$ is a positive integer. Considering equation (\ref{start})  modulo $t^6$, we find that 
\begin{equation} \label{loki}
1 \equiv \binom{m}{m-3}t^3 (1-t^3)^{(m-5)/3} + (1-t^3)^{(m-2)/3} \mod{t^6}.
\end{equation}
We have
$$
\binom{m}{m-3} = \binom{m}{3} = \binom{m}{2} \cdot \frac{m-2}{3} = (k(t^3-1)+1) \left( \frac{m-2}{3}  \right)
$$
and so
$$
1 \equiv (1-k) \left( \frac{m-2}{3} \right) t^3 + 1 - \left( \frac{m-2}{3} \right)  t^3 \mod{t^6}.
$$
It follows that
\begin{equation} \label{hammer}
k  \left( \frac{m-2}{3} \right)\equiv 0 \mod{t^3}.
\end{equation}
We claim that necessarily 
\begin{equation} \label{good}
m > \sqrt[3]{6} \, t^2.
\end{equation} 
Suppose not. Then
from $\binom{m}{2} =k(t^3-1)+1$,
$$
k  \left( \frac{m-2}{3} \right)  = \frac{ \binom{m}{2} -1}{t^3-1} \left( \frac{m-2}{3} \right) 
\leq  \frac{ \left( \sqrt[3]{6} \, t^2 (\sqrt[3]{6} \, t^2-1)-2 \right) \left( \sqrt[3]{6} \, t^2 -2 \right) }{6(t^3-1)}<  t^3,
$$
contradicting (\ref{hammer}) and the assumption that $m > 2$.

Combining this new bound  (\ref{good}) with (\ref{upperm}) now implies that $t \le 352$. A short computer search finds all possible pairs $(t, m)$ satisfying \eqref{upperm}, \eqref{smallhammer}, \eqref{hammer}, (\ref{good}) and $2 \leq t \leq 352$. The only $t$ with admissible $m$ are
\[
2 \le t \le 13 \quad\text{and}\quad t = 15,16,18,19,21,22,25
\]
where we find a number of pairs as follows.
$$
\begin{array}{|cc|cc|cc|cc|} \hline
t & \# \mbox{ of } m & t & \# \mbox{ of } m & t & \# \mbox{ of } m & t & \# \mbox{ of } m \\ \hline
2 & 6715 & 7 & 92 & 12 & 5 & 19 & 1 \\ 
3 & 3857 & 8 & 18 & 13 & 3 & 21 & 1 \\
4 & 1510 & 9 & 51 & 15 & 2 & 22 & 1 \\
5 & 59 & 10 & 39 & 16 & 9 & 25 & 1 \\ 
6 & 291 & 11 & 4 & 18 & 2 &  &  \\ \hline
\end{array}
$$
Checking the coefficient of $s^2$ in the expansion for $(t-s)^m$ for each of these pairs  $(t,m)$, we find no solutions to equation \eqref{masser} for $m>2$. Therefore \emph{all} solutions to equation \eqref{cubicnorm} are of the ``functional'' kind, given by
\[
(x,y) = (t^2,2t) \; \mbox{ and } \;  (t^2,-t).
\]
This completes our proof.



\end{document}